\documentclass{amsart}
\usepackage{amssymb}

\usepackage{color}

\newtheorem{theorem}{Theorem}

\newtheorem{lemma}{Lemma}

\newcommand{\pp}{\noindent {\em Proof. }}
\newcommand{\bee}[1]{\begin{equation}\label{#1}}
\newcommand{\beq}[1]{\begin{eqnarray}\label{#1}}
\newcommand{\ene}{\end{equation}}
\newcommand{\eqe}{\end{eqnarray}}

\begin{document}
\title{Graded identities of some simple Lie superalgebras}

\author[D. Repov\v s, and M. Zaicev]
{Du\v san Repov\v s and Mikhail Zaicev}

\address{Du\v san Repov\v s \\Faculty of Education, and
Faculty of  Mathematics and Physics, University of Ljubljana,
P.~O.~B. 2964, Ljubljana, 1001, Slovenia}
\email{dusan.repovs@guest.arnes.si}

\address{Mikhail Zaicev \\Department of Algebra\\ Faculty of Mathematics and
Mechanics\\  Moscow State University \\ Moscow,119992, Russia}
\email{zaicevmv@mail.ru}

\thanks{The first author was supported by the Slovenian Research Agency
grants P1-0292-0101 and J1-4144-0101. The second author was partially supported 
by RFBR, grant 13-01-00234a. We thank the referee for comments and suggestions.}



\begin{abstract}
We study $\mathbb{Z}_2$-graded identities of Lie superalgebras of the type
$b(t), t\ge 2$, over a field of characteristic zero. Our main result is that the
$n$-th codimension is strictly less than $(\dim b(t))^n$ asymptotically.
As a consequence we obtain an upper bound for ordinary (non-graded)
PI-exponent for each simple Lie superalgebra $b(t), t\ge 3$.
\end{abstract}


\maketitle

\section{Introduction}

In this paper we study numerical invariants of identities of Lie superalgebras.
One of the main numerical characteristics of the identities of an algebra $A$ over a 
field $F$ of characteristic zero is the sequence of codimensions $\{c_n(A)\}, n=1,2,
\ldots~~$, and its asymptotic behaviour. Many deep and interesting results in this 
area were proved during the last few  decades (see, for example, 
\cite{GZbook})  both in the
associative and the non-associative cases. In particular, in many classes of algebras
(associative \cite{GZ1}, Lie \cite{GRZ1}, \cite{GRZ2}, \cite{z2002}, Jordan, 
alternative and some others \cite{GShZ}) it was proved that if $A$ is a finite dimensional 
algebra, $\dim A=d$, and $F$ is algebraically closed then PI-exponent $exp(A)$ is 
equal to $d$ if and only if $A$ is simple. In general $exp(A)\le\dim A$ as it was
observed in \cite{BaDr},
\cite{GZ2}. Recently (see \cite{GZLond}) it was shown that $exp(L)<\dim L$ provided 
that $L$ is a simple Lie superalgebra of the type $b(t), t\ge 3$, (we use notations from
\cite{Sch} for simple Lie superalgebras). Unfortunately, an upper bound $\alpha = 
\alpha(t)<\dim b(t)$ was not found for $exp(b(t))$  in \cite{GZLond}.

Since any Lie superalgebra $L$ is $\mathbb{Z}_2$-graded, one can consider graded
codimensions $c_n^{gr}(L)$ and graded PI-exponent $exp^{gr}(L)$. Graded codimensions and 
graded PI-exponents of Lie superalgebras were studied earlier in several papers 
(see, for example, \cite{ZM1}, \cite{ZM2}, \cite{ZM3}). Existence and integrality 
of graded exponents were proved for some classes  of Lie superalgebras. On the other 
hand, there are no known examples where $exp^{gr}(L)$ is fractional.

There are some relations between graded and non-graded identities, codimensions and 
PI-exponents. In particular, 
\begin{equation}\label{ast}
c_n(A)\le c_n^{gr}(A)
\end{equation}
(see  \cite{Gia-Reg} or \cite{BaDr}) for any finite dimensional $G$-graded algebra
$A$ where $G$ is a finite group. Hence when $A$ is finite dimensional and simple and 
$exp(A)=\dim A$ it follows from (\ref{ast}) and \cite{BaDr} that $exp^{gr}(A)$ exists
and is equal to $\dim A$.

First series of examples with $exp(A)\ne\dim A$ where $A$ is a finite dimensional simple 
algebra is given by simple Lie superalgebras $b(t), t\ge 3$, of the dimension $\dim b(t)=2t^2-1$ 
(\cite{GZLond}). It is important to study asymptotics of $c_n^{gr}(b(t))$ and to compare 
it with the asymptotics of $c_n(b(t))$. The main result of this paper says that the 
(upper) graded PI-exponent of $b(t)$ is less than or equal to $t^2-1+t\sqrt{t^2-1}$. As a 
consequence of this result and (\ref{ast}) we obtain an upper bound for ordinary 
PI-exponent of $b(t)$, $exp(b(t))\le t^2-1+t\sqrt{t^2-1}$. In particular, the difference 
$\dim b(t) - exp(b(t))$ is at leasi $t^2-t\sqrt{t^2-1}$ which is a decreasing function of $t$
with limit $\frac{1}{2}$.

\section{Preliminaries}

Let $A$ be an algebra over a field $F$ of characteristic zero. Recall that $A$ is said 
to be $\mathbb{Z}_2$-graded algebra if $A$ has a vector space decomposition
$A=A_0\oplus A_1$ such that $A_0A_0+A_1A_1\subseteq A_0$, $A_0A_1+A_1A_0\subseteq A_1$.
Usually elements of $A_0$ are called even while elements of $A_1$ are called odd. Any
element of $A_0\cup A_1$ is called homogeneous. In particular, a Lie superalgebra $L$ 
is a $\mathbb{Z}_2$-graded algebra $L=L_0\oplus L_1$ satisfying the following two relations
$$
xy-(-1)^{|x||y|}yx=0,
$$
$$
x(yz)=(xy)z+(-1)^{|x||y|}y(xz)
$$
where $x,y,z$ are homogeneous elements and $|x|=0$ if $x$ is even while $|x|=1$ if $x$ 
is odd.

Denote by $\mathcal{L}(X,Y)$ a free Lie superalgebra with infinite sets of even
generators $X$ and odd generators $Y$. A polynomial $f=f(x_1,\ldots, x_m,y_1,\ldots,y_n) 
\in \mathcal{L}(X,Y)$ is said to be a graded identity of Lie superalgebra 
$L=L_0\oplus L_1$ if $f(a_1,\ldots, a_m,b_1,\ldots,b_n)=0$ whenever $a_1,\ldots, a_m\in 
L_0, b_1,\ldots,b_n\in L_1$.

Given positive integers $0\le k\le n$, denote by $P_{k,n-k}$ the subspace of all multilinear 
polynomials  $f=f(x_1,\ldots, x_k,y_1,\ldots,y_{n-k})\in \mathcal{L}(X,Y)$
of degree $k$ in even variables and of degree $n-k$ in odd variables. Denote by
$Id^{gr}(L)$ an ideal of $\mathcal{L}(X,Y)$ of all graded identities of $L$. Then
$P_{k,n-k}\cap Id^{gr}(L)$ is the subspace of all multilinear graded identities of $L$
of total degree $n$ depending on $k$ even variables and $n-k$ odd variables. Denote also
by $P_{k,n-k}(L)$ the quotient
$$
P_{k,n-k}(L)=\frac{P_{k,n-k}}{P_{k,n-k}\cap Id^{gr}(L)}.
$$
Then the graded $(k,n-k)$-codimension of $L$ is
$$
c_{k,n-k}(L)=\dim P_{k,n-k}(L)
$$
and the total graded codimension of $L$ is
\begin{equation}\label{f0}
c_n^{gr}(L)=\sum_{k=0}^n{n\choose k} c_{k,n-k}(L).
\end{equation}
If the sequence $\{c_n^{gr}(L)\}_{n\ge 1}$ is  exponentially bounded then one can
consider the related bounded sequence $\sqrt[n]{c_n^{gr}(L)}$. The latter sequence
has the following lower and upper limits 
$$
\underline{exp}^{gr}(L)=\liminf_{n\to\infty} \sqrt[n]{c_n^{gr}(L)},\qquad
\overline{exp}^{gr}(L)=\limsup_{n\to\infty} \sqrt[n]{c_n^{gr}(L)}
$$
called the lower  and upper PI-exponents of $L$, respectively. If an ordinary
limit exists, it is called an (ordinary) graded PI-exponent of $L$,
$$
exp^{gr}(L)=\lim_{n\to\infty} \sqrt[n]{c_n^{gr}(L)}
$$

Symmetric groups and their representations play an important role in the theory of 
codimensions. In particular, in the case of graded identities one can consider the 
$S_k\times S_{n-k}$-action on multilinear graded polynomials. Namely, the subspace
$P_{k,n-k}\subseteq\mathcal{L}(X,Y)$ has a natural structure of $S_k\times S_{n-k}$-
module where $S_k$ acts on even variables $x_1,\ldots,x_k$ while $S_{n-k}$ acts
on odd variables $y_1,\ldots,y_{n-k}$. Clearly, $P_{k,n-k}\cap Id^{gr}(L)$ is the
submodule under this action and we get an induced $S_k\times S_{n-k}$-action on 
$P_{k,n-k}(L)$. The character $\chi_{k,n-k}(L)=\chi(P_{k,n-k}(L))$ is called 
$(k,n-k)$ cocharacter of $L$. By Maschke's Theorem this character can be decomposed 
into the sum of irreducible characters
\begin{equation}\label{f1}
\chi_{k,n-k}(L)=\sum_{{\lambda\vdash k\atop \mu\vdash n-k}}
m_{\lambda,\mu}\chi_{\lambda,\mu}
\end{equation}
where $\lambda$ and $\mu$ are partitions of $k$ and $n-k$, respectively (all details
concerning representations of symmetric groups  can be found in \cite{JK}).
 
Recall that an irreducible $S_k\times S_{n-k}$-module with the character 
$\chi_{\lambda,\mu}$ is the tensor product of $S_k$-module with the character
$\chi_\lambda$ and $S_{n-k}$-module with the character $\chi_\mu$. In particular,
the dimension $\deg\chi_{\lambda,\mu}$ of this module is the product 
$d_\lambda d_\mu$ where $d_\lambda=\deg\chi_\lambda, d_\mu=\deg\chi_\mu$.
Taking into account multiplicities $m_{\lambda,\mu}$ in (\ref{f1}) we get the relation
\begin{equation}\label{f2}
c_{k,n-k}(L)=\sum_{{\lambda\vdash k\atop \mu\vdash n-k}}
m_{\lambda,\mu}d_\lambda d_\mu.
\end{equation}
A number of irreducible components in the decomposition of $\chi_{k,n-k}(L)$, i.e. 
the sum
$$
l_{k,n-k}(L)=\sum_{{\lambda\vdash k\atop \mu\vdash n-k}} m_{\lambda,\mu}
$$ 
is called the $(k,n-k)$-colength of $L$. If $\dim L<\infty$ then by Ado Theorem
(see \cite[Theorem 1.4.1]{Sch}), $L$ has a faithful finite dimensional graded
representation. Hence $L$ has an embedding $L\subset A=A_0\oplus A_1$ as a Lie
superalgebra where $A$ is a finite dimensional associative superalgebra. Given
$0\le k\le n$, consider the graded $(k,n-k)$-cocharacter of $A$:
$$
\chi_{k,n-k}(A)=\sum_{{\lambda\vdash k\atop \mu\vdash n-k}}
\overline m_{\lambda,\mu} \chi_{\lambda,\mu}.
$$
Then by \cite{Ber},
$$
\sum_{k=0}^n\sum_{{\lambda\vdash k\atop \mu\vdash n-k}}
\overline m_{\lambda,\mu} \le q(n)
$$
for some polynomial $q(n)$. Following the argument of the proof of 
\cite[Lemma 3.2]{GRZ1} we obtain that
$$
m_{\lambda,\mu} \le \overline m_{\lambda,\mu}.
$$
Hence in the finite dimensional case the total colength is polynomially bounded, 
that is for any $L$, $\dim L<\infty$, there exists a polynomial $f(n)$ such that
$$
\sum_{k=0}^n l_{k,n-k}(L) \le f(n).
$$
It follows that
\begin{equation}\label{f3}
c_{k,n-k}(L) \le f(n) d_\lambda^{max} d_\mu^{max}
\end{equation}
where $d_\lambda^{max}, d_\mu^{max}$ are maximal possible dimensions of $S_k$- and
$S_{n-k}$-representations, respectively, such that $m_{\lambda,\mu}\ne 0$. We
will use relation (\ref{f3}) for finding an upper bound for $\overline{exp}^{gr}(L)$.

\section{Dimensions of some $S_m$-representations}

In this section we prove some technical results which we will use later. Fix an integer
$t\ge 2$ and consider an irreducible $S_m$-representation with the character
$\chi_\mu, \mu=(\mu_1,\ldots,\mu_d), d\le t^2$. For convenience we will write
$\mu=(\mu_1,\ldots,\mu_{t^2})$ even in case $d<t^2$ assuming $\mu_{d+1}=\ldots =
\mu_{t^2}=0$.

We define the following function of a partition $\mu\vdash m$
\begin{equation}\label{f4}
\Phi(\mu)=\frac{1}{\left(\frac{\mu_1}{m}\right)^\frac{\mu_1}{m}\cdots 
\left(\frac{\mu_{t^2}}{m}\right)^\frac{\mu_{t^2}}{m}}
\end{equation}
In (\ref{f4}) we assume that $0^0=1$ if some of $\mu_j$ are equal to zero. The value of
$\Phi(\mu)^m$ is equal to $d_\mu$ up to a polynomial factor. More precisely, we have the 
following relation.
\begin{lemma}\label{l1}\cite[Lemma 1]{GZLond}
Let $m\ge 100$. Then
$$
\frac{\Phi(\mu)^m}{m^{t^4+t^2}} \le d_\mu \le m \Phi(\mu)^m.
$$
\end{lemma}
\hfill $\Box$

Now let $\lambda$ and $\mu$ be two partitions of $m$ with the corresponding Young diagrams
$D_\lambda, D_\mu$. We say that $D_\mu$ is obtained from $D_\lambda$ by pushing down one
box if there exist $1\le i<j\le t^2$ such that $\mu_i=\lambda_i-1, \mu_j=\lambda_j+1$
and $\mu_k=\lambda_k$ for all remaining $k$.

\begin{lemma}\label{l2}\rm{(see} \cite[Lemma 3]{GZLond}, \cite[Lemma 2]{ZaRe})
Let $D_\mu$ be obtained from $D_\lambda$ by pushing down one box. Then 
$\Phi(\mu) \ge \Phi(\lambda)$.
\end{lemma}
\hfill $\Box$

Now we define the weight of  partition $\mu=(\mu_1,\ldots, \mu_{t^2})$ as follows:
$$
wt~\mu=-(\mu_1+\cdots+\mu_{\frac{t^2-t}{2}}) + 
(\mu_{\frac{t^2-t}{2}+1}+\cdots+\mu_{t^2}).
$$

Recall (see\cite{GZbook}) that the hook partition $h(d,l,k)$ is a
partition with the Young diagram of the shape

\hskip2cm
\hskip2cm
\setlength{\unitlength}{2565sp}%

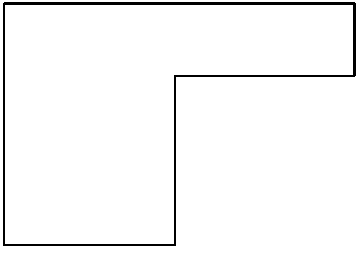
\begingroup\makeatletter\ifx\SetFigFont\undefined%
\gdef\SetFigFont#1#2#3#4#5{%
  \reset@font\fontsize{#1}{#2pt}%
  \fontfamily{#3}\fontseries{#4}\fontshape{#5}%
  \selectfont}%
\fi\endgroup%
\begin{picture}(4865,2677)(464,-1986)
\thinlines
{\color[rgb]{0,0,0}
\put(476,364){\line( 1, 0){3385}}
\put(3864,-336){\line(-1, 0){1738}}
\put(2126,-1486){\line( 0,-1){488}}
\put(2126,-1974){\line(-1, 0){1650}}
\put(476,-1974){\line( 0, 1){2338}}
}%
{\color[rgb]{0,0,0}
\put(2126,-344){\line( 0,-1){1463}}
}%
{\color[rgb]{0,0,0}\put(3864,-336){\line( 0, 1){700}}
}%
\put(5314, 14){\makebox(0,0)[lb]{\smash{{\SetFigFont{8}{9.6}{\rmdefault}{\mddefault}{\updefault}
}}}}
\put(3951,-949){\makebox(0,0)[lb]{\smash{{\SetFigFont{8}{9.6}{\rmdefault}{\mddefault}
{\updefault}
}}}}
\put(1226,451){\makebox(0,0)[lb]{\smash{{\SetFigFont{8}{9.6}{\rmdefault}{\mddefault}{\updefault}{\color[rgb]{0,0,0}
}%
}}}}
\put(2951,476){\makebox(0,0)[lb]{\smash{{\SetFigFont{8}{9.6}{\rmdefault}{\mddefault}{\updefault}{\color[rgb]{0,0,0}}%
}}}}
\put(4464,426){\makebox(0,0)[lb]{\smash{{\SetFigFont{8}{9.6}{\rmdefault}{\mddefault}{\updefault}{\color[rgb]{0,0,0}}%
}}}}
\put(2176,-1811){\makebox(0,0)[lb]{\smash{{\SetFigFont{8}{9.6}{\rmdefault}{\mddefault}{\updefault}{\color[rgb]{0,0,0}}%
}}}
}
\end{picture}%
\vskip .2cm
Here the first $d$ rows have length $l+k$ and remaining $k$ rows have length $l$. 
We slightly modify this notion and say that a partition $\mu=(\mu_1,\ldots,\mu_{t^2})\vdash m$ is a hook $h(s,r)$ if $\mu_1=\ldots=\mu_{\frac{t^2-t}{2}}=s$ and
$\mu_{\frac{t^2-t}{2}+1}=\cdots=\mu_{t^2}=r<s$.

The following observation is elementary.
\begin{lemma}\label{l3}
Let $m$ be a multiple of $t(t^2-1)$. Then there exists a hook partition $\mu=h(s,r)$
of $m$ with $s=r\frac{t+1}{t-1}$ and $wt~\mu=0$.
\end{lemma}
\pp Let $m=it(t^2-1)$. If we take $\mu=h(r,s)$ with $s=(t+1)i$, $r=(t-1)i$ then 
the number of boxes in the first $\frac{t^2-t}{2}$ rows, that is
$\mu_1+\cdots+\mu_{\frac{t^2-t}{2}}$, equals to
$$
s\frac{t^2-t}{2}=it\frac{(t-1)(t+1)}{2}=\frac{m}{2}.
$$
Similarly, the number of boxes in all remaining rows of $D_\mu$ equals to
$$
r\frac{t^2+t}{2}=it\frac{(t-1)(t+1)}{2}=\frac{m}{2}.
$$
Hence $wt~\mu=0$ and we are done.
\hfill $\Box$

\begin{lemma}\label{l4}
Let $m$ be a multiple of $t(t^2-1)$ and let $\mu=h(s,r)$ be the hook partition
with zero weight as in Lemma \ref{l3}. Then $\Phi(\mu)=t\sqrt{t^2-1}$.
\end{lemma}
\pp  Since
$$
\frac{\mu_1}{m}=\ldots=\frac{\mu_{\frac{t^2-t}{2}}}{m}=\frac{s}{m},\qquad
\frac{\mu_{\frac{t^2-t}{2}+1}}{m}=\ldots=\frac{\mu_{t^2}}{m}=\frac{r}{m}
$$
and $m=rt(t+1), s=r\frac{t+1}{t-1}$, we have
$$
\frac{r}{m}=\frac{1}{t(t+1)},\qquad  \frac{s}{m}=\frac{1}{t(t-1)}.
$$
Hence
$$
\Phi(\mu)=\frac{1}{\left(\frac{1}{t(t+1)}\right)^\frac{t^2+t}{2t(t+1)}
\left(\frac{1}{t(t-1)}\right)^\frac{t^2-t}{2t(t-1)}}=
\left(t^2(t+1)(t-1)\right)^\frac{1}{2}=t\sqrt{t^2-1}.
$$
\hfill $\Box$

For an arbitrary partition of weight zero we have the following.

\begin{lemma}\label{l5}
Let $m$ be a multiple of $t(t^2-1)$ and let $\nu$ be a partition of $m$
with $wt~\nu=0$. Then $\Phi(\nu)\le t\sqrt{t^2-1}$.
\end{lemma}
\pp The Young diagram $D_\nu$  of $\nu$ consists of two parts. The first one $\overline\nu$
contains first $\frac{t^2-t}{2}$ rows and the second part $\overline{\overline{\nu}}$
contains all remaining rows. Pushing down boxes inside $\overline\nu$ and 
$\overline{\overline{\nu}}$ separately we get new partition $\nu'\vdash m$ with
$wt~\nu'=0$ maximally close to hook partition. That is, first $0<i\le \frac{t^2-t}{2}$
rows of $D_{\nu'}$ have the length $a$ and rows $i+1,\ldots,\frac{t^2-t}{2}$ (in case
$i< \frac{t^2-t}{2}$) have the length $a-1$. Similarly,
$$
\nu'_{\frac{t^2-t}{2}+1}=\ldots=\nu'_{\frac{t^2-t}{2}+j}=b,\quad
\nu'_{\frac{t^2-t}{2}+j+1}=\ldots=\nu'_{t^2}=b-1
$$
for some $j$. But under our assumption $m$ admits a hook partition by Lemma \ref{l3},
hence $\frac{m}{2}$ is a multiple of $\frac{t^2-t}{2}$. It follows that 
$i=\frac{t^2-t}{2}$. Similarly, $j=\frac{t^2+t}{2}$ and $\nu'= h(a,b)$. Finally
note that if $m$ admits a hook partition $\mu=h(r,s)$  of weight zero then $\mu$ is
uniquely defined. Hence $a=b\frac{t+1}{t-1}$ and $\Phi(\nu')=t\sqrt{t^2-1}$ by Lemma
\ref{l4}. By applying Lemma \ref{l2} we complete the proof.
\hfill $\Box$

The main goal of this section is to get a similar upper bound for $\Phi(\mu)$
for any $\mu\vdash m$ without any restriction on $m$ and with $wt(\mu) \le 1$.

First we prove an easy technical result.

\begin{lemma}\label{d1}
Let $\lambda=(\lambda_1,\ldots,\lambda_s)$ be a partition of $n$ such that
$\lambda_1-\lambda_s\ge 2s$. Then by pushing down one or more boxes in $D_\lambda$
one can get a partition $\mu=(\mu_1,\ldots \mu_s)\vdash n$ with $\mu_s=\lambda_s+1$ and $\mu_1>\lambda_1-s$.
Similarly, one can get $\nu=(\nu_1,\ldots \nu_s)\vdash n$ with $\nu_1=\lambda_1-1$ and $\nu_s<\lambda_s+s$.
\end{lemma}
 
\pp First we find $\mu$. If $s=2$ then the statement is obvious. Suppose $s>2$. 
Then we push down boxes in $D_\lambda$ using only rows $2,3,\ldots, s$. If we get on some step the diagram $D_\mu$
with $\mu_s=\lambda_s+1$ then we proof is completed. Otherwise we will get a diagram $D_{\bar \mu}$ where
$\bar \mu_1=\lambda_1, \bar \mu_2=\cdots=\bar \mu_t=p+1$, $\bar \mu_{t+1}=\cdots=\bar \mu_s=p$ for some $p$
and some $2\le t\le s$. Moreover, $p=\lambda_s$ if $t<s$ or $p+1=\lambda_s$ if $t=s$. In this case we can 
cut $s-1$ boxes from the first row of  $D_{\bar \mu}$ and the glue one box to each row $2,\ldots, s$ in $D_{\bar \mu}$.
Then the partition $\mu=(\mu_1,\ldots,\mu_s)$, $\mu_1=\bar\mu_1-s+1, \mu_j=\bar\mu_j+1, j=2,\ldots,s$, 
satisfies all conditions and we are done.

Similarly, if we push down boxes only in rows $1,\ldots, s-1$ in $D_\lambda$ then either we will get a partition
$\nu=(\nu_1,\ldots,\nu_s)$ with $\nu_1=\lambda_1-1,\nu_s=\lambda_s$ on some step or we will get a partition
$\bar\nu=(\bar\nu_1,\ldots,\bar\nu_s)$ such that$\bar\nu_1=\cdots=\bar\nu_t=p+1$, $\bar\nu_{t+1}=\cdots=\bar\nu_{s-1}=p$,
$\bar\nu_s=\lambda_s$ for some $1\le t\le s-1$. In the latter case we push down one box from each row $1,\ldots, s-1$
to the last row of $D_{\bar\nu}$. Then we get the required $\nu\vdash n$ and the proof is completed.

\hfill $\Box$

Now we consider partitions with $t^2$ components whose weight cannot be increased by 
pushing down boxes in the Young diagram.

\begin{lemma}\label{d2}
Let $\mu=(\mu_1,\ldots,\mu_{t^2})$ be a partition whose weight cannot be increased 
by pushing down boxes. Then $\mu_1-\mu_{t^2} \le 4t^2$ and $wt(\mu)\ge -2t^4$.
\end{lemma}

\pp Denote $p=\frac{t^2+t}{2}, q=\frac{t^2-t}{2}$ for brevity. Clearly, $\mu_q\le 
\mu_{q+1}+1$. If $\mu_1-\mu_q\ge 3q$ then by pushing down
boxes we can get a partition $\mu'=(\mu_1',\ldots,\mu_{q}')$ with $\mu_q'=
\mu_q+2$ by Lemma \ref{d1}. Hence $\mu_1-\mu_q< 3q$. Similarly, we can
get $(\mu_{q+1}'',\ldots,\mu_{q+p}'')$ with $\mu_{q+1}''=\mu_{q+1}-2$ provided
that $\mu_{q+1}-\mu_{q+p} \ge 3p$. Therefore $\mu_{q+1}-\mu_{q+p}<3p$.
Finally we obtain
$$
\mu_1-\mu_{q+p}<3p+3q+1=3t^2+1< 4t^2.
$$

For proving the second part of our lemma we split $D_\mu$ into two parts $D_1$ and $D_2$
where $D_1$ consists of the first $q$ rows of $D_\mu$ while $D_2$ consists of the last $p$ 
rows of $D_\mu$. By our assumption we cannot cut one box from $D_1$ and glue it
to $D_2$. Denote by $a$ and $b$ the number of boxes in $D_1,D_2$, respectively.
Denote also $\mu_{p+q}=x$. By the first part of the lemma $\mu_1\le 4t^2+x$.
Hence $a\le (4t^2+x)q$. Obviously, $b\ge px$. Hence
$$
wt(\mu)=b-a\ge x(p-q)-4t^2q\ge -4t^2\frac{t^2-t}{2}\ge -2t^4
$$
and we complete the proof.

\hfill $\Box$

Next lemma shows how to reduce this problem to the case $wt~\mu=0$ and $m=jt(t^2-1)$.

\begin{lemma}\label{l6}
Let $\mu=(\mu_1,\ldots,\mu_{t^2})$ be a partition of $m$ and let $wt~\mu\le 1$.
Then there exist an integer $m_0\ge m$ and a partition $\nu\vdash m_0$ such that
\begin{itemize}
\item[1)] $m_0-m\le 6t^6$,

\item[2)]   $wt~\nu=0$,
\item[3)]   $m_0$ is a multiple of $t^2(t-1)$,
\item[4)]   
$$
\Phi(\mu)\le (m+6t^6)^{(\frac{t^4+t^2+2}{m})^{6t^6}}\Phi(\nu).
$$
\end{itemize}
\end{lemma}
\pp 
First we reduce the question to the case $wt~\mu=0$. If $wt~\mu=1$ then we can
add one extra box to the first row of $D_\mu$ and get a partition of zero weight. 

Let $wt(\mu)<0$. By Lemma \ref{l2} and Lemma \ref{d2} we can suppose that 
$wt(\mu) \ge -2t^4$. If we add one box to each of rows $1,2,\ldots,t^2-t+1$ of $D_\mu$
we get the Young diagram $D_{\rho}$ of partition $\rho\vdash m+t^2-t+1$ with  
$wt(\rho)=wt(\mu)+1$. Applying this procedure at most $2t^4$ times we get 
$\rho'\vdash m_0'$ with $wt(\rho')=0$ where
$$
m_0'\le m+(t^2-t+1)\cdot 2t^4 \le m+4t^6.
$$

If $m_0'$ is a multiple of $t(t^2-1)$ then there is nothing to do.
Otherwise there exists $0<i<t(t^2-1)$ such that $m_0'+i$ is a multiple of $t(t^2-1)$.
Note that $m_0'$ is even since it admits a partition of weight zero. Hence $i$  is 
also even.

First we enlarge $D_{\rho'}$ to $D_{\mu'}$ by adding $\frac{t^2-1}{2}$ boxes to all
$t^2-t$ first rows. Then also $wt~\mu'=0$. Since $\mu'_{t^2-t}-\mu'_{t^2-t+1} \ge 
\frac{t^2-1}{2}$, we can glue $\frac{i}{2}<t\frac{t^2-1}{2}$ boxes to the last $t$
rows of $D_{\mu'}$ and get $D_{\mu''}$. Finally, we glue $\frac{i}{2}$ boxes to the 
first row of $D_{\mu''}$ and obtain the diagram $D_\nu$ such that $wt~\nu=0$. 
Denote by $m_0$ the number of boxes of $D_\nu$.
As follows from our procedure, an upper bound for $m_0$ is
$$
m+4t^6+(t^2-t)\frac{t^2-1}{2}+t(t^2-1)<6t^6+m.
$$

It is shown in \cite[Lemma 7]{ZaRe} that if $\lambda\vdash n-1, \lambda = 
(\lambda_1,\ldots, \lambda_{d})$, $\lambda'\vdash n, \lambda' = 
(\lambda_1',\ldots, \lambda_{d}')$ and $D_\lambda$ is obtained from  $D_{\lambda'}$ 
by cutting one box then
$$
\Phi(\lambda)\le n^\frac{d^2+d+2}{n}\Phi(\lambda').
$$ 
Hence
$$
\Phi(\mu)\le (m+6t^6)^{(\frac{t^4+t^2+2}{m})^{6t^6}}\Phi(\nu)
$$
and we complete the proof.
\hfill $\Box$

As a corollary of Lemma \ref{l5} and Lemma \ref{l6} we immediately obtain

\begin{lemma}\label{l7}
Let $\mu=(\mu_1,\ldots,\mu_{t^2})$ be a partition of $m$ and let $ wt~\mu\le 1$.
Then there exists a polynomial $g(m)$ such that $\Phi(\mu)\le g(m)^\frac{1}{m} t\sqrt{t^2-1}$.
\end{lemma}

\hfill $\Box$

\section{Graded codimensions of Lie superalgebras of type {\it b(t})}

In this section we use notations from \cite{Sch}. Recall that $L=b(t), t\ge 2$, is
a Lie superalgebra of $2t\times 2t$ matrices of the type
$$
\left(
  \begin{array}{cc}
    A & B \\
    C & -A^T \\
  \end{array}
\right),
$$
where $A,B,C\in M_t(F)$, $B^T=B, C^T=-C$ and $tr A=0$.
Here  the map $X\to X^T$ is the transpose involution. Decomposition $L=L_0\oplus L_1$
is defined by setting
$$
L_0 = \left\lbrace \left(
           \begin{array}{cc}
             A & 0 \\
             0 & -A^T \\
           \end{array}
         \right)
 \mid A\in M_t(F), tr(A)=0  \right\rbrace,
$$
and
$$
L_1 = \left\lbrace \left(
           \begin{array}{cc}
             0 & B \\
             C & 0 \\
           \end{array}
         \right)
 \mid B^T=B, C^T=-C\in M_t(F)  \right\rbrace.
$$
Super-Lie product on $L$ is given by
$$
[x,y]=xy-(-1)^{|x||y|}yx
$$
for homogeneous $x,y\in L_0\cup L_1$. 

It is not difficult to see that also $L$ has $\mathbb{Z}$-grading
\begin{equation}\label{f6}
L=L^{(-1)}\oplus L^{(0)}\oplus L^{(1)}
\end{equation}
where $L^{(0)}= L_0$, 
\begin{equation}\label{f7}
L^{(-1)}= \left\lbrace \left(
           \begin{array}{cc}
             0 & 0 \\
             C & 0 \\
           \end{array}
         \right)
 \mid C^T=-C\in M_t(F)  \right\rbrace,
\end{equation}
\begin{equation}\label{f8}
L^{(1)} = \left\lbrace \left(
           \begin{array}{cc}
             0 & B \\
             0 & 0 \\
           \end{array}
         \right)
 \mid B^T=B\in M_t(F)  \right\rbrace
\end{equation}
and $L^{(n)}=0$ for all $n\ne0,\pm 1$. In particular, $L^{(-1)}\oplus L^{(1)}=L_1$
and $\dim L^{(0)}=t^2-1,\dim L^{(-1)}=\frac{t(t-1)}{2},\dim L^{(1)}=\frac{t(t+1)}{2}$.

Let $\chi_{k.n-k}(L)$ be $(k,n-k)$-cocharacter of $L$. Consider its decomposition
(\ref{f1}) into irreducible components.

\begin{lemma}\label{l8}
Let $m_{\lambda,\mu}\ne 0$ in (\ref{f1}). Then $D_\lambda$ lies in the strip of width $t^2-1$, that is,
 $\lambda=(\lambda_1,\ldots,\lambda_d)$ with $d\le t^2-1$. In particular, $d_\lambda\le  \alpha(k) (t^2-1)^k$ 
for some polynomial $\alpha(k)$ .
\end{lemma}
\pp
Denote $A=FS_k$. Recall that, given a partition $\lambda=(\lambda_1,\ldots,\lambda_d) \vdash k$,
the irreducible $S_k$-module corresponding to $\lambda$ is isomorphic to the minimal
left ideal generated by an essential idempotent $e_{T_\lambda}$ constructed in the following way.

Let $T_\lambda$ be Young tableau that is Young diagram $D_\lambda$ filled up by integers 
$1,\ldots,k$. Denote by $R_{T_\lambda}$ and $C_{T_\lambda}$ row and column stabilizers
in $S_k$ of $T_\lambda$, respectively. Then
$$
R(T_\lambda)=\sum_{\sigma\in R_{T_\lambda}} \sigma~, \quad
C(T_\lambda)=\sum_{\tau\in C_{T_\lambda}} ({\rm sgn}\tau)\tau
$$
and
$$
e_{T_\lambda}= R(T_\lambda) C(T_\lambda).
$$
It is known that $e_{T_\lambda}^2=\alpha e_{T_\lambda}$, $0\ne \alpha\in \mathbb{Q}$,
and an irreducible $FS_k$-module $M$ has the character $\chi_\lambda$ if and only if
$e_{T_\lambda}M\ne 0$. In particular, if $M$ is an irreducible 
$FS_k\times FS_{n-k}$-submodule in $P_{k,n-k}(L)$ with the character $\chi_{\lambda,\mu}$
then $M$ can be generated by a multilinear polynomial of the type
$e_{T_\lambda}\varphi(x_1,\ldots,x_k,y_1,\ldots,y_{n-k})$ with even $x_1,\ldots,x_k$ 
and odd $y_1,\ldots,y_{n-k}$ (since $M$ is the direct sum of isomorphic irreducible 
$S_k$-modules with characters $\chi_\lambda$). 
From the relation $e_{T_\lambda}^2=
\alpha e_{T_\lambda}\ne 0$ it follows that the polynomial
$$
\psi(x_1,\ldots,x_k,y_1,\ldots,y_{n-k})=C(T_\lambda)
e_{T_\lambda}\varphi(x_1,\ldots,x_k,y_1,\ldots,y_{n-k})
$$
also generates $M$.

Suppose now that $d>t^2-1$. Then $D_\lambda$ contains at least one column of height $d$
greater than $t^2-1=\dim L_0$. In this case $\psi$ depends on at least one alternating set
of even variables of order greater than $\dim L_0$. Standard arguments show that in this 
case $\psi$ is an identity of $L$, a contradiction. Hence $d\le t^2-1$. Now by
\cite[Lemma 6.2.5]{GZbook} there exists a polynomial $\alpha(k)$ such that
$d_\lambda\le\alpha(k)(t^2-1)^k$  and  we complete the proof.
\hfill $\Box$

\begin{lemma}\label{l9}
Let $m_{\lambda,\mu}\ne 0$ in (\ref{f1}). Then $ wt~\mu\le 1$.
\end{lemma}
\pp
As in the previous lemma an irreducible $FS_k\times FS_{n-k}$-submodule $M$ of
$P_{k,n-k}(L)$ with the character $\chi_{\lambda,\mu}$ can be generated by
$$
\psi=\psi(x_1,\ldots,x_k,y_1,\ldots,y_{n-k})=C(T_\mu)e_{T_\mu}\varphi
$$
for some multilinear polynomial $\varphi$. The set of variables $\{y_1,\ldots,y_{n-k}\}$
can be split into disjoint union
$$
\{y_1,\ldots,y_{n-k}\}=Y_1\cup\ldots\cup Y_p
$$
where $p=\mu_1$, every set $Y_j$ consists of odd indeterminates with the indices from
the $j$-th column of $T_\mu$. In particular, $\psi$ is alternating on any subset $Y_j$, 
$1\le j \le p$, and we cannot substitute the same basis elements of $L$ instead of
distinct variables from the same column of $T_\mu$, otherwise the value of $\psi$
will be zero. Hence the minimal degree in $\mathbb{Z}$-grading (\ref{f6}),(\ref{f7}), 
(\ref{f8}) of the value of $\psi$ on $L$ is equal to $q=wt~\mu$. So, if $q>1$ then $\psi$  
is an identity of $L$ since $L^{(q)}\oplus L^{(q+1)}\oplus\cdots=0$, a contradiction.

\hfill $\Box$

Now we are ready to prove the main result of the  paper.

\begin{theorem}\label{t1}
Let $L$ be a Lie superalgebra of the type $b(t), t\ge 2$, over a field $F$ of 
characteristic zero. Then there exists a polynomial $h=h(n)$ such that
$$
c_n^{gr}(L)\le h(n)(t^2-1+t\sqrt{t^2-1})^n.
$$
In particular,
$$
\overline{exp}^{gr}(L) \le t^2-1+t\sqrt{t^2-1}< 2t^2-1=\dim L.
$$
\end{theorem}

\pp Consider the inequality (\ref{f3}) for $c_{k,n-k}(L)$. By Lemma \ref{l8}, 
$d_\lambda^{max}\le\alpha(k) (t^2-1)^k$ and by Lemma \ref{l9} we have $wt~\mu\le 1$ where
$d_\mu=d_\mu^{max}$. Then by Lemma \ref{l1} and Lemma \ref{l7},
$$
d_\mu^{max} \le(n-k)g(n-k) (t\sqrt{t^2-1})^{n-k}.
$$
Hence
$$
c_{k,n-k}(L)\le f(n) (n-k)\alpha(k)g(n-k)(t^2-1)^k (t\sqrt{t^2-1})^{n-k}.
$$
Clearly one can take a polynomial $h'=h'(n)$ such that $\alpha(k) g(n-k)\le h'(n)$ for all
$k=0,\ldots, n$. Then
$$
c_{k,n-k}(L)\le h(n) (t^2-1)^k (t\sqrt{t^2-1})^{n-k}
$$
where $h(n)=nf(n)h'(n)$. Now by (\ref{f0})
$$
c_n^{gr}(L)\le h(n)\sum_{k=0}^n {n \choose k} (t^2-1)^k (t\sqrt{t^2-1})^{n-k}=
h(n)(t^2-1+t\sqrt{t^2-1})^n.
$$

Obviously,
$$
\overline{exp}^{gr}(L)=\limsup_{n\to\infty}\sqrt[n]{c_n^{gr}(L)}
\le t^2-1+t\sqrt{t^2-1}
$$
and we complete the proof of Theorem \ref{t1}.
\hfill $\Box$

As a consequence of Theorem \ref{t1} we get an upper bound for ordinary PI-exponent
of $L=b(t)$, $t \ge 2$.

\begin{theorem}\label{t2}
Let $L$ be a  Lie superalgebra of the type $L=b(t)$, $t \ge 2$, over a field of 
characteristic zero. Then $exp(L)\le t^2-1+t\sqrt{t^2-1}$.
\end{theorem}
\pp The statement easily follows from Theorem \ref{t1}, the inequality 
$c_n^{gr} \le c_n(L)$ (\cite{BaDr}, \cite{Gia-Reg}) and from the existence of $exp(L)$
(\cite{GZLond}).
\hfill $\Box$

\end{document}